\magnification 1100
\def\arXiv{oui}
  
 \def\oui{oui} 
  
\ifx\arXiv\oui
\else
 \pdfpagewidth=210truemm
 \pdfpageheight=297truemm 
\fi
  
  %
  %

  %
  %

  \catcode`@=12 

 \def\defrefnote#1{\definexref{#1}{{\the\footnotenumber}}{refnotes}}

  %
  %


\ifx\couleurs\oui
\input graphicx
 \pdfpagewidth=210truemm
 \pdfpageheight=297truemm 
 \voffset=-5mm
\fi

\input eplain.tex
\expandafter\def\expandafter\newdimen\expandafter{\newdimen}

\ifx\couleurs\oui
\beginpackages
\usepackage{color}
\endpackages 
 \pdfpagewidth=210truemm
 \pdfpageheight=297truemm 
\long\def\rge#1{{\color{red}#1}}

\definecolor{bleu-iecn}{cmyk}{.98,.13,.1,.55}

\else
\long\def\rge#1{#1}

\fi

\makeatletter
\def\numberedfootnote{%
ÊÊ\global\advance\footnotenumber by 1
ÊÊ\@eplainfootnote{{\number\footnotenumber}}%
}%
\def\makecolumns#1/#2 {\par \begingroup
ÊÊ \@columndepth = #1
ÊÊ \advance\@columndepth by -1
ÊÊ \divide \@columndepth by #2
ÊÊ \advance\@columndepth by 1
ÊÊ \@linestogoincolumn = \@columndepth
ÊÊ \@linestogo = #1
ÊÊ \currentcolumn = 1
ÊÊ \def\@endcolumnactions{%
ÊÊÊÊÊÊ\ifnum \@linestogo<2
ÊÊÊÊÊÊÊÊ \the\crtok \egroup \endgroup \par 
ÊÊÊÊÊÊ\else
ÊÊÊÊÊÊÊÊ \global\advance\@linestogo by -1
ÊÊÊÊÊÊÊÊ \ifnum\@linestogoincolumn<2
ÊÊÊÊÊÊÊÊÊÊÊÊ\global\advance\currentcolumn by 1
ÊÊÊÊÊÊÊÊÊÊÊÊ\global\@linestogoincolumn = \@columndepth
ÊÊÊÊÊÊÊÊÊÊÊÊ\the\crtok
ÊÊÊÊÊÊÊÊ \else
ÊÊÊÊÊÊÊÊÊÊÊÊ&\global\advance\@linestogoincolumn by -1
ÊÊÊÊÊÊÊÊ \fi
ÊÊÊÊÊÊ\fi
ÊÊ }%
ÊÊ \makeactive\^^M
ÊÊ \letreturn \@endcolumnactions
ÊÊ \@columnwidth = \hsize
ÊÊÊÊ \advance\@columnwidth by -\parindent
ÊÊÊÊ \divide\@columnwidth by #2
ÊÊ \penalty\abovecolumnspenalty
ÊÊ \noindent 
ÊÊ \valign\bgroup
ÊÊÊÊ &\hbox to \@columnwidth{\strut \hsize = \@columnwidth ##\hfil}\cr
}%
\makeatother

\lefteqnumbers
   \def\testd{oui}
   \def\choixlat{\ifx\numadroite\testd\righteqnumbers
            \else  \lefteqnumbers\fi}
    \choixlat

\catcode`@=\letter
\def\@eplainfootnote#1{\let\@sf\empty 
  \ifhmode\edef\@sf{\spacefactor\the\spacefactor}\/\fi
  \global\advance\hlfootlabelnumber by 1
  \hlstart@impl{foot}{\hlfootlabel}%
  \hldest@impl{footback}{\hlfootbacklabel}%
  \hbox{$^{(#1)}$}%
  \hlend@impl{foot}%
  \@sf\vfootnote{#1.}%
}%
\catcode`@=\other

  \interfootnoteskip=0pt
  \let\note=\numberedfootnote
  \everyfootnote={\eightpoint\leftskip=5truemm\rightskip5truemm}
  
  \hsize150truemm\vsize 240truemm\hoffset=5truemm

  \pretolerance=500\tolerance=1000\brokenpenalty=5000
  \parindent3mm
  
  \countdef\temps=170
  \temps=\time
  \countdef\nminutes=171{\nminutes = \time}
  \countdef\nheure=172
  \def\heure{\begingroup                     
     \temps = \time \divide\temps by 60
     \nheure = \temps                        
     \nminutes = \time
     \multiply\temps by 60
     \advance\nminutes by -\temps            
     \ifnum\nminutes<10 \toks1 = {0}%
     \else\toks1 = {}%
     \fi
     \number\nheure h\the\toks1 \number\nminutes  
  \endgroup}%

  \newcount\chstart
  \chstart=\pageno
 \headline={\ifnum\pageno=\chstart {\hfill} \else {\hss \tenrm --\ \folio\ --\hss}\fi}
  \footline={\hfill}
  \normalbaselines
  \frenchspacing
    \def\dater{\vglue-10mm\rightline{(\the\day/\the\month/\the\year)}}
  \def\dateheure{\vglue-10mm\rightline{(\the\day/\the\month/\the\year,\ \heure)}}

  \newif\ifpagetitre \pagetitretrue
\newtoks\hautpagetitre \hautpagetitre={\hfill}
\newtoks\baspagetitre \baspagetitre={\hfill}
\newtoks\auteurcourant \auteurcourant={\hfill}
\newtoks\titrecourant \titrecourant={\hfill}
\newtoks\hautpagegauche
\newtoks\hautpagedroite
\newtoks\hautpagemilieu
\hautpagemilieu={\tenrm\hfil -- \folio\ -- \hfil}
\hautpagegauche={\ifx\midfolio\oui\the\hautpagemilieu\else\tenrm\folio\hfill\the\auteurcourant\hfill\fi}
\hautpagedroite={\ifx\midfolio\oui\the\hautpagemilieu\else\hfill\the\titrecourant\hfill\tenrm\folio\fi}
\newtoks\baspagegauche \baspagegauche={\hfil}
\newtoks\baspagedroite \baspagedroite={\hfil}
\headline={\ifpagetitre\the\hautpagetitre
\else\ifodd\pageno\the\hautpagedroite\else\the\hautpagegauche\fi\fi }
\footline={\ifpagetitre\the\baspagetitre
\else\ifodd\pageno\the\baspagedroite
\else\the\baspagegauche\fi\fi \global\pagetitrefalse}

\def\pageblanche{\vfill\eject\pagetitretrue
\null\vfill\eject
\pagetitretrue
}
\def\chgtpage{\ifodd\pageno \else
\pageblanche \fi \pagetitretrue\titreun=0\footnotenumber=0}

\def\chgtpageincrtitreun{\ifodd\pageno \else
\pageblanche \fi \pagetitretrue\footnotenumber=0}

\def\majnombres{\ifodd\pageno \else
\pageblanche \fi \pagetitretrue\hautpoly\titreun=0\footnotenumber=0}

\def\hautspages#1#2{\auteurcourant={\ninepcap#1}\titrecourant={\nineit#2}}

\ifnum\chstart=\pageno \pagetitretrue\fi
  

  \def\PAR{\par}
  
  \def\leftnote#1{\vadjust{\setbox1=\vtop{\hsize 20mm\parindent=0pt\eightpoint
  \baselineskip=9pt\rightskip=4mm plus 4mm\vskip-4.7mm#1}\hbox{\kern-2cm\smash{\box1}}}}

  
  \def\raggedcenter{\leftskip=20pt plus 10em  
       \rightskip=\leftskip 
        \parfillskip=0pt 
         \spaceskip=.3333em \xspaceskip=.5em 
          \pretolerance=9999 \tolerance=9999
           \hyphenpenalty=9999 \exhyphenpenalty=9999 }
           
  \def\titrecentre#1{{\parindent0mm\raggedcenter
       \spaceskip=.6em plus .2em minus .2em\xspaceskip=.6em plus .2em minus .2em
        \tit#1\par}}
        


  \def\oui{oui}
  
\def\fontetitreun{\ifx\paradouze\oui\douzepts\gpdouze\twelvebf\textfont1=\twelveib\else
\quatorzepts\gpquatorze\fourteenbf\fi}

\def\fontetitreunl{\douzepts\textfont1=\twelveib\scriptfont1=\tenib\fourteenti}
 
 \def\fontetitredeux{\textfont1=\eleveni\ifx\paradouze\oui\onzepts\scriptfont1=\ninei\elevenit\else
                        \douzepts\twelveit\fi}
 
   \def\fontetitredeuxb{\ifx\paradouze\oui\onzepts\eleventi\gponze\textfont1=\elevenib\scriptfont1=\nineib
                         \else\douzepts\twelveti\scriptfont1=\twelveib\scriptfont1=\tenib\gpdouze\fi}
                         
\def\fontetitredeuxl{\onzepts\textfont1=\elevenbf\scriptfont1=\ninebf\twelvebf}
  
\def\fontetitretrois{\textfont0=\elevenrm\scriptfont0=\eightrm\textfont1=\eleveni
                      \scriptfont1=\eighti\scriptscriptfont1=\sixi\elevenit}
                      
\def\fontetitrequatre{\textfont0=\elevenrm\scriptfont0=\eightrm\textfont1=\eleveni
                      \scriptfont1=\eighti\scriptscriptfont1=\sixi\elevenrm}
  
  \newcount\titreun\titreun=0
  \newcount\titredeux\titredeux=0
  \newcount\titretrois\titretrois=0
  \newcount\titrequatre\titrequatre=0
  \newcount\enonce\enonce=0
  
  \def\incr#1{\global\advance#1 by 1 {\the #1}}
  \def\avance#1{\global\advance#1 by 1}
  \def\init#1{\global#1=0}
  
  \long\def\Indentation#1#2{\setbox10=\hbox{\fontetitreun#1}
  	                    \ifdim\wd10 < 4mm
                         \setbox10=\hbox to 4mm{\box10\hfill}
                       \else\ifdim\wd10 < 6mm
                         \setbox10=\hbox to 6mm{\box10\hfill}
  	                    \else\ifdim\wd10 < 8mm
                         \setbox10=\hbox to 8mm{\box10\hfill}
                       \else\ifdim\wd10 < 12mm
                         \setbox10=\hbox to 12mm{\box10\hfill}
                       \fi\fi\fi\fi
                       \dimen10=\hsize
                       \advance \dimen10 by -\wd10
                       \noindent \box10 %
                       \ignorespaces
                       \hbox{\vtop{\hsize=\dimen10\raggedright\noindent\fontetitreun#2}}}

  \long\def\paraun#1{\removelastskip\par\medskip\goodbreak\vskip0pt plus.01\vsize\penalty-100
                \vskip0pt plus-.01\vsize
  	              \init{\titredeux}\ifnum\optionparag=1{\init\eqnumber\init\enonce}\else{}\fi
                  \goodbreak{\fontetitreun
  	                \Indentation{\incr{\titreun}.\ }{\fontetitreun #1\par}}\nobreak\medskip}

 %
 %
 \long\def\paraunc#1{\removelastskip\par\bigskip\goodbreak\vskip0pt plus.01\vsize\penalty-100
                \vskip0pt plus-.01\vsize
  	              \init{\titredeux}
                 \ifnum\optionparag=1{\init{\eqnumber}\init\enonce}\else{}\fi
                  \goodbreak
  	                {\parindent0mm\raggedcenter\fontetitreun\incr{\titreun}.\ 
                     \fontetitreun #1\par}\nobreak\medskip}
                     
\newtoks\titreunl
\titreunl={\ifnum\titreun=1{I}\fi%
\ifnum\titreun=2{II}\fi%
\ifnum\titreun=3{III}\fi%
\ifnum\titreun=4{IV}\fi%
\ifnum\titreun=5{V}\fi%
\ifnum\titreun=6{VI}\fi%
\ifnum\titreun=7{VII}\fi%
\ifnum\titreun=8{VIII}\fi%
\ifnum\titreun=9{IX}\fi%
\ifnum\titreun=10{X}\fi%
\ifnum\titreun=11{XI}\fi%
\ifnum\titreun=12{XII}\fi%
\ifnum\titreun=13{XIII}\fi%
}
\long\def\paraunl#1{\removelastskip\par\bigskip\bigskip\goodbreak\vskip0pt plus.01\vsize\penalty-100
                \vskip0pt plus-.01\vsize
  	              \init{\titredeux}\ifnum\optionparag=1{\init\eqnumber\init\enonce}\else{}\fi
                  \goodbreak{\fontetitreunl
  	                \Indentation{\global\advance\titreun by 1{\the\titreunl}.\ }{\fontetitreunl #1\par}}\nobreak\smallskip}

  
  \long\def\paradeux#1{\init{\titretrois}\vskip0pt plus.01\vsize\penalty-10
                \vskip0pt plus-.01\vsize\ifx \elie\oui\medskip\ifnum\titredeux>0\medskip\fi\fi
                 \Indentation{\fontetitredeux\the\titreun${\cdot}$\incr{\titredeux}.}
                              {\fontetitredeux\textfont1=\eleveni#1}\nobreak\par }
  
  \long\def\paradeuxb#1{\init{\titretrois}\vskip0pt plus.001\vsize\penalty-10
                \vskip0pt plus-.01\vsize{\ifx \elie\oui\medskip\ifnum\titredeux>0\medskip\fi\fi
                  \Indentation
  {\fontetitredeuxb\the\titreun${\cdot}$\incr{\titredeux}.}{ \fontetitredeuxb#1}}\nobreak
\smallskip}

\newtoks\titredeuxl
\titredeuxl={\ifnum\titredeux=1{A}\fi%
\ifnum\titredeux=2{B}\fi%
\ifnum\titredeux=3{C}\fi%
\ifnum\titredeux=4{D}\fi%
\ifnum\titredeux=5{E}\fi%
\ifnum\titredeux=6{F}\fi%
\ifnum\titredeux=7{G}\fi%
\ifnum\titredeux=8{H}\fi%
\ifnum\titredeux=9{I}\fi%
\ifnum\titredeux=10{J}\fi%
\ifnum\titredeux=11{K}\fi%
\ifnum\titredeux=12{L}\fi%
\ifnum\titredeux=13{M}\fi%
}
 \long\def\paradeuxl#1{\init{\titretrois}\vskip0pt plus.001\vsize\penalty-10
                \vskip0pt plus-.01
                \vsize \bigskip%
                  \Indentation
     {\fontetitredeuxl\global\advance\titredeux by 1
  \quad \the\titreunl${\cdot}$\the\titredeuxl.}{ \fontetitredeuxl#1}
  \removelastskip\nobreak\smallskip}
  

  \long\def\paratrois#1{\init{\titrequatre}\ifdim\lastskip<\smallskipamount
                \removelastskip\smallskip\fi
                 \vskip0pt plus.01\vsize\penalty-10
                  \vskip0pt
plus-.01\vsize{\ifx \elie\oui\ifnum\titretrois>0\medskip\fi\fi
\Indentation{\fontetitretrois\the\titreun${\cdot}$\the\titredeux${\cdot}$\incr{\titretrois}.\ }
  {\hskip0mm\baselineskip=14pt\fontetitretrois#1}\nobreak\smallskip}}
  
  
  \long\def\paratroisl#1{\init{\titrequatre}\ifdim\lastskip<\smallskipamount
                \removelastskip\fi
                 \vskip0pt plus.01\vsize\penalty-10
                  \vskip0pt
plus-.01\vsize\ifx \elie\oui\bigskip
\fi
\Indentation{\fontetitretrois\quad \quad \the\titreunl{${\cdot}$}\the\titredeuxl${\cdot}$\incr{\titretrois}.\ }
  {\hskip0mm\fontetitretrois#1}\nobreak\smallskip}


  \long\def\paraquatre#1{\ifdim\lastskip<\smallskipamount
                \removelastskip\smallskip\fi
                 \vskip0pt plus.01\vsize\penalty-10
                  \vskip0pt
                  plus-.01\vsize\par
 
\Indentation{\fontetitrequatre \the\titreun{${\cdot}$}\the\titredeux${\cdot}$\the\titretrois${\cdot}$\incr{\titrequatre}.\ }
{\hskip0mm\fontetitrequatre#1}\nobreak\smallskip}


\newtoks\titrequatrel
\titrequatrel={\ifnum\titrequatre=1{a}\fi%
\ifnum\titrequatre=2{b}\fi%
\ifnum\titrequatre=3{c}\fi%
\ifnum\titrequatre=4{d}\fi%
\ifnum\titrequatre=5{e}\fi%
\ifnum\titrequatre=6{f}\fi%
\ifnum\titrequatre=7{g}\fi%
\ifnum\titrequatre=8{h}\fi%
\ifnum\titrequatre=9{i}\fi%
\ifnum\titrequatre=10{j}\fi%
\ifnum\titrequatre=11{k}\fi%
\ifnum\titrequatre=12{l}\fi%
\ifnum\titrequatre=13{m}\fi%
}
\long\def\paraquatrel#1{\ifdim\lastskip<\smallskipamount
                \removelastskip\smallskip\fi
                 \vskip0pt plus.01\vsize\penalty-10
                  \vskip0pt
                  plus-.01\vsize{\bigskip
\Indentation{\global\advance\titrequatre by 1
\fontetitrequatre\quad \quad \quad \the\titreunl${\cdot}$\the\titredeuxl${\cdot}$\the\titretrois${\cdot}$\the\titrequatrel.\ }
{\hskip0mm\fontetitrequatre#1}\nobreak\smallskip}}

\ifx\optionkeys\oui
\def\drefun#1{\definexref{¤#1}{{\the\titreun}}{}} 
\def\drefdeux#1{\definexref{¤#1}{{\the\titreun}.{\the\titredeux}}{}}
\def\dreftrois#1{\definexref{¤#1}{{\the\titreun}.{\the\titredeux}.{\the\titretrois}}{}}
\else
\def\drefun#1{\definexref{prg#1}{{\the\titreun}}{}} 
\def\drefdeux#1{\definexref{prg#1}{{\the\titreun}.{\the\titredeux}}{}}
\def\dreftrois#1{\definexref{prg#1}{{\the\titreun}.{\the\titredeux}.{\the\titretrois}}{}}
\fi

%


  \long\def\propdeux#1#2#3#4{%
       \avance{\enonce}
       \leavevmode\edef\temp{#2}%
         \ifx\temp\empty 
          \else
           \definexref{#2}{#1~{\the\titreun.\the\enonce}}{enonces}
            \definexref{s#2}{{\the\titreun.\the\enonce}}{enonces}
             \fi
\smallskip
      \noindent{\bf#1\ {\bf\the\titreun.\the\enonce{#3}.}\enspace}{\sl#4\par}%
      \ifdim\lastskip<\medskipamount \removelastskip\penalty55\par \fi
   }

  \long\def\propun#1#2#3#4{%
      \avance{\enonce}
       \leavevmode\edef\temp{#2}%
        \ifx\temp\empty 
          \else
           \definexref{#2}{#1~{\the\enonce}}{enonces}
            \definexref{{s#2}}{{\the\enonce}}{enonces}
             \fi
   \par 
     \noindent{\bf#1\ {\bf\the\enonce{#3}.}\enspace}{\sl#4\par}%
     \ifdim\lastskip<\medskipamount \removelastskip\penalty55\medskip\fi
  }
  
  \long\def\prop#1#2#3#4{\ifnum\optionparag=1
                          \propdeux{#1}{#2}{\textfont1=\elevenib#3}{#4} \else\propun{#1}{#2}{\textfont1=\elevenib#3}{#4}\fi}

  \long\def\propt#1#2#3{\ifx\tpf\oui \prop{Th\'eo\-r\`eme}{#1}{#2}{#3}\par
                       \else\prop{Theorem}{#1}{#2}{#3}\par\fi}
  \long\def\Propt#1#2{\propt{#1}{}{#2}}
  \long\def\propl#1#2#3{\ifx\tpf\oui\prop{Lem\-me}{#1}{#2}{#3}\par
                         \else\prop{Lemma}{#1}{#2}{#3}\par\fi}
  \long\def\Propl#1#2{\propl{#1}{}{#2}}
  \long\def\propc#1#2#3{\ifx\tpf\oui\prop{Corol\-laire}{#1}{#2}{#3}\par
                         \else\prop{Corollary}{#1}{#2}{#3}\par\fi}
  
  \long\def\propp#1#2#3{\prop{Pro\-po\-si\-tion}{#1}{#2}{#3}\par}
  \long\def\Propp#1#2{\propp{#1}{}{#2}} 
  \long\def\propd#1#2#3{\ifx\tpf\oui\prop{D\'efi\-nition}{#1}{#2}{#3}\par
                       \else\prop{Definition}{#1}{#2}{#3}\par\fi} 
  
  \long\def\proptd#1#2#3{\ifx\tpf\oui\prop{Th\'eor\`eme et d\'efi\-nition}{#1}{#2}{#3}\par
                       \else\prop{Theorem and definition}{#1}{#2}{#3}\par\fi}


  
  \newcount\optionparag\optionparag=1
  
  \long\def\section#1#2{\ifnum\optionparag=1 \paraun{#2} 
                        \else\goodbreak{\fontetitreun
  	                \Indentation{#1.\ }{#2}}\nobreak\smallskip\fi}

  \def\eqconstruct#1{\ifnum\optionparag=1{\the\titreun\hbox{$\cdot$}#1}\else{#1}\fi}

  
  
  \def\numref{oui}  
  
  \newcount\mesref\mesref=0 
  \def\defbib#1{\ifx\numref\oui\global\advance\mesref by 1 \definexref{#1}{{\the
                 \mesref}}{}\else\definexref{#1}{#1}{}\fi}
  \def\bibtem#1{\defbib{#1}\item{\citer{#1}}}
  \def\citer#1{[\ref{#1}]}
  \def\citeplus#1#2{[\ref{#1}; #2]}

  
  \font\seventeenmsa=msam10 at 17pt    
  \font\fourteenmsa=msam10 at 14pt
  \font\twelvemsa=msam10 at 12pt
  \font\tenmsa=msam10                 
  \font\ninemsa=msam10 at 9pt 
  \font\eightmsa=msam10 at 8pt 
  \font\sevenmsa=msam7 
  \font\sixmsa=msam10 at 6pt
  \font\fivemsa=msam5
  \newfam\msafam\textfont\msafam=\tenmsa\scriptfont\msafam=\sevenmsa\scriptscriptfont\msafam=\fivemsa
  
  \font\seventeenbb=msbm10 at 17pt     
  \font\fourteenbb=msbm10 at 14pt
  \font\twelvebb=msbm10 at 12pt
  \font\tenbb=msbm10                   
  \font\ninebb=msbm10 at 9pt 
  \font\eightbb=msbm10 at 8pt 
  \font\sevenbb=msbm7 
  \font\sixbb=msbm10 at 6pt
  \font\fivebb=msbm5 
  \newfam\bbfam\textfont\bbfam=\tenbb\scriptfont\bbfam=\sevenbb\scriptscriptfont\bbfam=\fivebb
  \def\bb{\fam\bbfam\tenbb}%

  \font\seventeenscaln=eusm10 at 17pt   
  \font\twelvescaln=eusm10 at 12pt
  \font\tenscaln=eusm10                
  \font\ninescaln=eusm10 scaled 900
  \font\eightscaln=eusm10 scaled 800
  \font\sevenscaln=eusm10 scaled 700
  \font\sixscaln=eusm10 scaled 600
   
  \newfam\scalnfam\textfont\scalnfam=\tenscaln\scriptfont\scalnfam=\sevenscaln\scriptscriptfont\scalnfam=\sixscaln
  \def\scaln{\fam\scalnfam\tenscaln}%
  \def\scal{\scaln}
  
  \font\tenscalb=eusb10                

  \font\sevenscalb=eusb10 scaled 700

  \newfam\scalbfam\textfont\scalbfam=\tenscalb\scriptfont\scalbfam=\sevenscalb
  %
  
  %
  %
  \font\fourteenrm=cmr12 scaled 1200
  \font\elevenrm=cmr10 at 11pt
  \font\twelverm=cmr12
  \font\ninerm=cmr9
  \font\eightrm=cmr8      
  \font\sevenrm=cmr7
  \font\sixrm=cmr6

  \font\seventeenpcap=cmcsc10 at 17pt
  \font\tenpcap=cmcsc10                        
  \font\ninepcap=cmcsc9
  \font\eightpcap=cmcsc8
  \font\sevenpcap=cmcsc10 scaled 700
  
  \newfam\pcapfam\textfont\pcapfam=\tenpcap\scriptfont\pcapfam=\sevenpcap
  \def\pcap{\fam\pcapfam\tenpcap}
  
  \font\seventeenrm=cmbx12 scaled 1400

  \font\fourteenbf=cmbx10 scaled 1400
  
  \font\twelvebf=cmbx12
  \font\elevenbf=cmbx10 at 11pt
  \font\ninebf=cmbx9  
  \font\eightbf=cmbx8
  \font\sixbf=cmbx6
  
  \font\tengot=eufm10                           
   
  \font\eightgot=eufm10 at 8truept 
  \font\sevengot=eufm7 
  \font\sixgot=eufm10 at 6 truept 
   
  \newfam\gotfam
  \textfont\gotfam=\tengot\scriptfont\gotfam=\sevengot\scriptscriptfont\gotfam=\sixgot
  \def\got{\fam\gotfam\tengot}%

  
  \def\tit{%
  \textfont0=\seventeenrm\scriptfont0=\tenrm\def\rm{\fam0\seventeenrm}%
  \textfont1=\seventeenib\scriptfont1=\twelveib%
  \textfont2=\seventeensy\scriptfont2=\twelvesy\scriptscriptfont2=\ninesy
  \textfont3=\seventeenex
  \textfont\itfam=\seventeenti
  \def\it{\fam\itfam\seventeenti}%
  \textfont\bbfam=\seventeenbb \scriptfont\bbfam=\twelvebb
  \def\bb{\fam\bbfam\seventeenbb}%
  \textfont\msafam=\seventeenmsa\scriptfont\msafam=\twelvemsa
  \textfont\scalnfam=\seventeenscaln
  \def\pcap{\fam\pcapfam\seventeenpcap}
  \normalbaselineskip=25pt\normalbaselines\rm}

  \font\seventeenti=cmbxti10 scaled 1680
  
  \font\fourteenti=cmbxti10 at 14pt
  
  \font\twelveti=cmbxti10 scaled 1200
  \font\eleventi=cmbxti10 at 11pt

  %
  %
  \font\twelveit=cmti12	
  \font\elevenit=cmti10 scaled 1100
  \font\nineit=cmti9
  \font\eightit=cmti8
  \font\sevenit=cmti7

  %
  %
 
 \font\seventeenib=cmmib10 scaled 1680
  \font\fourteenib=cmmib10 scaled 1400
  \font\twelveib=cmmib10 scaled 1200
  \font\elevenib=cmmib10 scaled 1100
  \font\tenib=cmmib10
\font\eightib=cmmib10 scaled 800
  \font\nineib=cmmib10 scaled 900
\font\sevenib=cmmib10 scaled 700
\font\sixib=cmmib10 scaled 600
\font\fiveib=cmmib10 scaled 500

\ifx\ITAN\oui
\else
\innernewfam\cmmibfam
\textfont\cmmibfam=\tenib
\scriptfont\cmmibfam=\sevenib
\scriptscriptfont\cmmibfam=\fiveib
\def\ib{\fam\cmmibfam\tenib}
\fi

  %
  %
  \font\twelvei=cmmi10 scaled 1200
  \font\eleveni=cmmi10 scaled 1100
  \font\ninei=cmmi9
  \font\eighti=cmmi8 
  \font\seveni=cmmi7 			                
  \font\sixi=cmmi6
  
  \font\ninesl=cmsl9                    
  \font\eightsl=cmsl8 
  \font\sevensl=cmsl10 at 7pt

  \font\ninett=cmtt9                    
  \font\eighttt=cmtt8
  \font\seventt=cmtt10 scaled 700

  \font\seventeensy=cmsy10 scaled 1680    
  \font\fourteensy=cmsy10 scaled 1400
  \font\twelvesy=cmsy10 scaled 1176
  
  \font\ninesy=cmsy9                      
  \font\eightsy=cmsy8
  \font\sixsy=cmsy6
  \font\seventeenex=cmex10 at 17pt
  \font\fourteenex=cmex10 at 14pt
  \font\twelveex=cmex10 at 12pt
  \font\nineex=cmex10 at 9pt
  \font\eightex=cmex10 at 8pt
  \font\sevenex=cmex10 at 7pt
  \font\sixex=cmex10 at 6pt
  \font\fiveex=cmex10 at 5pt
  
   
  \font\fourteengp=cmmi10 at 14pt
  
  \font\twelvegp=cmmib10 at 12pt
  \font\elevengp=cmmib10 at 11pt
  \font\tengp=cmmib10                          
  \font\ninegp=cmmib10 at 9pt 
  \font\eightgp=cmmib8 
   
  \font\sixgp=cmmib6


  \def\gponze{\textfont0=\elevenbf\scriptfont0=\eightbf\scriptscriptfont0=\sixbf
           \textfont1=\elevengp\scriptfont1=\eightgp\scriptscriptfont1=\sixgp}
  \def\gpdouze{\textfont0=\twelvebf\scriptfont0=\tenbf\scriptscriptfont0=\ninebf
           \textfont1=\twelvegp\scriptfont1=\tengp\scriptscriptfont1=\ninegp}        
  
 \def\gpquatorze{\textfont0=\fourteenbf\scriptfont0=\twelvebf\scriptscriptfont0=\elevenbf
           \textfont1=\fourteengp\scriptfont1=\twelvegp\scriptscriptfont1=\elevengp}

  
  \expandafter\chardef\csname pre amssym.def at\endcsname=\the\catcode`\@
  \catcode`\@=11
  \def\undefine#1{\let#1\undefined}
  \def\newsymbol#1#2#3#4#5{\let\next@\relax
   \ifnum#2=\@ne\let\next@\msafam@\else
   \ifnum#2=\tw@\let\next@\bbfam@\fi\fi
   \mathchardef#1="#3\next@#4#5}
  \def\mathhexbox@#1#2#3{\relax
   \ifmmode\mathpalette{}{\m@th\mathchar"#1#2#3}%
   \else\leavevmode\hbox{$\m@th\mathchar"#1#2#3$}\fi}
  \def\hexnumber@#1{\ifcase#1 0\or 1\or 2\or 3\or 4\or 5\or 6\or 7\or 8\or
   9\or A\or B\or C\or D\or E\or F\fi}
  
  \def\setboxz@h{\setbox\z@\hbox}
  \def\wdz@{\wd\z@}
  \def\boxz@{\box\z@}
  
  \edef\msafam@{\hexnumber@\msafam}
  \mathchardef\dabar@"0\msafam@39
  
  \edef\bbfam@{\hexnumber@\bbfam}
  \def\widehat#1{\setboxz@h{$\m@th#1$}%
   \ifdim\wdz@>\tw@ em\mathaccent"0\bbfam@5B{#1}%
   \else\mathaccent"0362{#1}\fi}
  \def\widetilde#1{\setboxz@h{$\m@th#1$}%
   \ifdim\wdz@>\tw@ em\mathaccent"0\bbfam@5D{#1}%
   \else\mathaccent"0365{#1}\fi}
  \newsymbol\leqq 1335          
  \newsymbol\leqslant 1336
  \newsymbol\lessgtr 1337       
  \newsymbol\backprime 1038     
  \newsymbol\risingdotseq 133A  
  \newsymbol\fallingdotseq 133B 
  \newsymbol\succcurlyeq 133C   
  \newsymbol\geqq 133D          
  \newsymbol\geqslant 133E
  \newsymbol\nmid 232D
  \newsymbol\nexists 2040
  \newsymbol\smallsetminus 2272
  \newsymbol\varnothing 203F 
  \catcode`\@=\active

  \catcode`\@=11
  \newcount\typofr\typofr=1
  
  \catcode`\;=\active
  \def;{\ifnum\typofr=1\relax\ifhmode\ifdim\lastskip>\z@\unskip\fi
     \kern.2em\fi\string;\else\string;\fi}
  
  \catcode`\:=\active
  \def:{\ifnum\typofr=1\relax\ifhmode\ifdim\lastskip>\z@\unskip\fi
  \penalty\@M\ \fi\string:\else\string:\fi}
  
  \catcode`\!=\active
  \def!{\ifnum\typofr=1\relax\ifhmode\ifdim\lastskip>\z@\unskip\fi
     \kern.2em\fi\string!\else\string!\fi}
  
  \catcode`\?=\active
  \def?{\ifnum\typofr=1\relax\ifhmode\ifdim\lastskip>\z@\unskip\fi
     \kern.2em\fi\string?\else\string?\fi}

  \def\francais{\typofr=1\def\tpf{oui}}
  \def\anglais{\typofr=2\def\tpf{non}\def\english{oui}}
  \def\oui{oui}
  \francais
  
  \catcode`\@=12
  

%
\def\raye #1|{\leavevmode\setbox1=\hbox{#1}%
\raise .5pt\hbox to \wd1{\xleaders\hbox{\rge{$ \sct / $}%
\kern 1pt}\hfill\kern -1pt }\kern -\wd1 \unhbox1\relax }

\def\barre #1|{\leavevmode\setbox1=\hbox{#1}%
\rlap{\color{red}\vrule height 2.4pt depth -1.2pt width \wd1}\color{black} \unhbox1\relax}

  

  
  \def\og{\leavevmode\raise.24ex\hbox{$\scriptscriptstyle\langle\!\langle\>$}}    
  \def\fg{\leavevmode\raise.24ex\hbox{$\scriptscriptstyle\>\rangle\!\rangle$}}    

  \def\d{\,{\rm d}}
  \def\dd{{\rm d}}

  \def\CC{{\bb C}}
  \def\N{{\bb N}}

  \def\PP{{\bb P}}

  \def\HH{{\scal H}}

  \def\O{{\scal O}}
  \def\P{{\scaln P}}
  \def\QQ{{\scaln Q}}

  \def\frac#1#2{{#1\over #2}}
  \def\di#1#2{\sct#1\atop{\sct#2}}
  \def\tri#1#2#3{{\sct#1\atop\sct#2}\atop\sct#3}

  \def\qedbox{$\rlap{$\sqcap$}\sqcup$}           
  \def\qed{\nobreak\hfill\penalty250 \hbox{}\nobreak\hfill\qedbox\par }

  \def\pnt{prime number theorem}

  \def\¤{\S\thinspace}

  \def\¥{$\bullet$ }
  
  
  \def\e{{\rm e}}
  
  \def\md#1#2{\equiv#1\,({\rm mod\,}#2)}

  \def\epsilon{\varepsilon}

  \def\phi{\varphi}
  \def\theta{\vartheta}
  \def\rho{\varrho}
  \def\dm{{\textstyle{1\over 2}}}
  \def\txt{\textstyle}
  \def\dsp{\displaystyle}
  \def\sct{\scriptstyle}
  \def\pf{\noi{\it Proof. }}
  \def\nid{\ifnum\typofr=1\par\noindent{\it D\'emonstration. }\else\pf\fi}
  \def\noi{\noindent}
  \def\rem{\ifnum\typofr=1\noi{\it Remarque.}\ \else\noi{\it Remark.}\ \fi}
  \def\rems{\ifnum\typofr=1\noi{\it Remarques.}\ \else\noi{\it Remarks.}\ \fi}

  \def\sset{\smallsetminus}

  \def\1{{\bf 1}}
  \def\|{\Vert}

  \def\leq{\leqslant}
  \def\geq{\geqslant}
  \def\wh{\widehat}

  \def\eg{{e.g.}}
  \newsymbol\subsetneqq 2324
  \newsymbol\subsetneq 2328

  \def\log{\mathop{\rm log}\nolimits}
  \def\ft#1#2{{\txt{#1\over #2}}}


  
\def\vbs#1{\left|#1\right|}


  \def\pmb#1{\setbox0=\hbox{#1}%
  \kern-.025em\copy0\kern-\wd0\kern.05em\copy0\kern-\wd0\kern-.025em\raise .0433em\box0 }

  
  \skewchar\eighti='177 \skewchar\sixi='177
  \skewchar\eightsy='60 \skewchar\sixsy='60
  
  \def\eightpoint{%
  \textfont0=\eightrm\scriptfont0=\sixrm\scriptscriptfont0=\fiverm
  \def\rm{\fam0\eightrm}%
  \textfont1=\eighti\scriptfont1=\sixi
  \scriptscriptfont1=\fivei\def\oldstyle{\fam1\seveni}%
  \textfont2=\eightsy\scriptfont2=\sixsy\scriptscriptfont2=\fivesy
  \textfont3=\eightex\scriptfont3=\sixex
  \textfont\itfam=\eightit
  \def\it{\fam\itfam\eightit}%
  \textfont\slfam=\eightsl
  \def\sl{\fam\slfam\eightsl}%
  \textfont\bbfam=\eightbb \scriptfont\bbfam=\sixbb\scriptscriptfont\bbfam=\fivebb
  \def\bb{\fam\bbfam\eightbb}%
  \textfont\msafam=\eightmsa\scriptfont\msafam=\sixmsa
  \textfont\scalnfam=\eightscaln
  \def\scaln{\fam\scalnfam\eightscaln}
  \textfont\ttfam=\eighttt
  \def\tt{\fam\ttfam\eighttt}%
\textfont\gotfam=\eightgot
  \textfont\bffam=\eightbf\scriptfont\bffam=\sixbf\scriptscriptfont\bffam=\fivebf
  \def\bf{\fam\bffam\eightbf}%
  \ifx\ITAN\oui\else\textfont\cmmibfam=\eightib
       \scriptfont\cmmibfam=\sixib
        \scriptscriptfont\cmmibfam=\fiveib
         \def\ib{\fam\cmmibfam\eightib}
   \fi
  \textfont\pcapfam=\eightpcap
  \def\pcap{\fam\pcapfam\eightpcap}
  \abovedisplayskip=2pt plus2pt minus 2pt
  \belowdisplayskip=2pt plus1pt minus 2pt
  \abovedisplayshortskip= 1pt plus 2pt minus 1pt
  \belowdisplayshortskip= 1pt plus 2pt minus 1pt
  \smallskipamount=2pt plus 1pt minus 2pt
  \medskipamount=3pt plus 2pt minus 2pt
  \bigskipamount=7pt plus 3pt minus 3pt
  \setbox\strutbox=\hbox{\vrule height 5pt depth 2pt width 0pt}%
  \normalbaselineskip=9pt\normalbaselines\rm}

  \def\({\left(}
  \def\){\right)}
  
  \def\footnoterule{\kern -2pt\hrule width 7truecm\kern 2.4pt}
  
  \def\xnotedef#1{\definexref{#1}{\noexpand\number\footnotenumber}{Note}}%

  
  
  \def\ninepoint{%
  \textfont0=\ninerm\scriptfont0=\sixrm\scriptscriptfont0=\fiverm
  \def\rm{\fam0\ninerm}%
  \textfont1=\ninei\scriptfont1=\sixi
  \scriptscriptfont1=\fivei\def\oldstyle{\fam1\ninei}%
  \textfont2=\ninesy\scriptfont2=\sixsy\scriptscriptfont2=\fivesy
  \textfont3=\nineex\scriptfont3=\sixex
  \textfont\itfam=\nineit
  \def\it{\fam\itfam\nineit}%
  \textfont\slfam=\ninesl
  \def\sl{\fam\slfam\ninesl}%
  \textfont\bbfam=\ninebb\scriptfont\bbfam=\sixbb\scriptscriptfont\bbfam=\fivebb
  \def\bb{\fam\bbfam\ninebb}%
  \textfont\msafam=\ninemsa\scriptfont\msafam=\sixmsa\scriptscriptfont\msafam=\fivemsa
  \textfont\scalnfam=\ninescaln
  \def\scaln{\fam\scalnfam\ninescaln}
  \textfont\ttfam=\ninett
  \def\tt{\fam\ttfam\ninett}%
  \textfont\bffam=\ninebf\scriptfont\bffam=\sixbf\scriptscriptfont\bffam=\fivebf
  \def\bf{\fam\bffam\ninebf}%
  \abovedisplayskip=3pt plus2pt minus 2pt
  \belowdisplayskip=3pt plus1pt minus 2pt
  \abovedisplayshortskip= 2pt plus 2pt minus 1pt
  \belowdisplayshortskip= 2pt plus 2pt minus 1pt
  \smallskipamount=2pt plus 1pt minus 2pt
  \medskipamount=3pt plus 2pt minus 2pt
  \bigskipamount=7pt plus 3pt minus 3pt
  \setbox\strutbox=\hbox{\vrule height 5pt depth 2pt width 0pt}%
  \normalbaselineskip=11pt plus.3pt minus.2pt\normalbaselines\rm}

  \def\sevenpoint{%
  \textfont0=\sevenrm\scriptfont0=\sixrm\scriptscriptfont0=\fiverm
  \def\rm{\fam0\sevenrm}%
  \textfont1=\seveni\scriptfont1=\sixi
  \scriptscriptfont1=\fivei\def\oldstyle{\fam1\seveni}%
  \textfont2=\sevensy\scriptfont2=\sixsy\scriptscriptfont2=\fivesy
  \textfont3=\sevenex\scriptfont3=\fiveex
  \textfont\itfam=\sevenit
  \def\it{\fam\itfam\sevenit}%
  \textfont\slfam=\sevensl
  \def\sl{\fam\slfam\sevensl}%
  \textfont\bbfam=\sevenbb \scriptfont\bbfam=\sixbb\scriptscriptfont\bbfam=\fivebb
  \def\bb{\fam\bbfam\sevenbb}%
  \textfont\msafam=\sevenmsa\scriptfont\msafam=\sixmsa
  \textfont\scalnfam=\sevenscaln
  \def\scaln{\fam\scalnfam\sevenscaln}
  \textfont\bffam=\sevenbf\scriptfont\bffam=\sixbf\scriptscriptfont\bffam=\fivebf
  \def\bf{\fam\bffam\sevenbf}%
  \textfont\ttfam=\seventt
  \abovedisplayskip=2pt plus2pt minus 2pt
  \belowdisplayskip=2pt plus1pt minus 2pt
  \abovedisplayshortskip= 1pt plus 2pt minus 1pt
  \belowdisplayshortskip= 1pt plus 2pt minus 1pt
  \smallskipamount=2pt plus 1pt minus 2pt
  \medskipamount=3pt plus 2pt minus 2pt
  \bigskipamount=7pt plus 3pt minus 3pt
  \setbox\strutbox=\hbox{\vrule height 5pt depth 2pt width 0pt}%
  \normalbaselineskip=9pt\normalbaselines\rm}

 \def\onzepts{%
 \textfont0=\elevenrm\scriptfont0=\ninerm
 \textfont1=\eleveni\scriptfont1=\ninei
}

\def\douzepts{%
  \textfont0=\twelverm\scriptfont0=\tenrm\def\rm{\fam0\twelverm}%
  \textfont1=\twelvei\scriptfont1=\teni%
  \textfont2=\twelvesy\scriptfont2=\tensy\scriptscriptfont2=\eightsy
  \textfont3=\twelveex
  \textfont\itfam=\twelveti
  \def\it{\fam\itfam\twelveti}%
  \textfont\bffam=\twelvebf\scriptfont\bffam=\tenbf\scriptscriptfont\bffam=\eightbf
  \def\bf{\fam\bffam\twelvebf}%
  \textfont\bbfam=\twelvebb \scriptfont\bbfam=\tenbb
  \def\bb{\fam\bbfam\twelvebb}%
  \textfont\msafam=\twelvemsa\scriptfont\msafam=\tenmsa
  \textfont\scalnfam=\twelvescaln
  \normalbaselineskip=15pt\normalbaselines\rm}

\def\quatorzepts{%
  \textfont0=\fourteenrm\scriptfont0=\twelverm\def\rm{\fam0\fourteenrm}%
  \textfont1=\fourteenib\scriptfont1=\twelveib%
  \textfont2=\fourteensy\scriptfont2=\twelvesy\scriptscriptfont2=\tensy
  \textfont3=\fourteenex
  \textfont\itfam=\fourteenti
  \def\it{\fam\itfam\fourteenti}%
  \textfont\bffam=\fourteenbf\scriptfont\bffam=\twelvebf\scriptscriptfont\bffam=\tenbf
  \def\bf{\fam\bffam\fourteenbf}%
  \textfont\bbfam=\fourteenbb \scriptfont\bbfam=\twelvebb
  \def\bb{\fam\bbfam\fourteenbb}%
  \textfont\msafam=\fourteenmsa\scriptfont\msafam=\twelvemsa
  \textfont\scalnfam=\twelvescaln
  \normalbaselineskip=18pt\normalbaselines\rm}


\def\AA{{\it Acta Arith.}}

\def\picture #1 by #2 (#3){\leavevmode\vbox to #2{
     \hrule width #1 height 0pt depth 0pt
      \vfill
       \special{picture #3}}}

\def\illustration #1 by #2 (#3) scaled #4{\dimen1=#2
  \divide\dimen1 by 1000
  \multiply\dimen1 by #4
  \vtop to \dimen1{\dimen1=#1
  \divide\dimen1 by 1000
  \multiply\dimen1 by #4
  \hsize=\dimen1\vss
  \noindent\special{illustration #3 scaled #4}}}

\ifx\couleurs\oui

\fi

\optionparag=1
\vsize=255truemm
\voffset=-5truemm
\ifx\optionkeymacros\undefined\else \fi

\catcode`\Œ=\active\defŒ{{\aa}}       
\catcode`\º=\active\defº{\int}        
\catcode`\=\active\def{\c c}        
\catcode`\¶=\active\def¶{\partial}    
\catcode`\Ä=\active\defÄ{\oint}       
\catcode`\Æ=\active\defÆ{\triangle}   
\catcode`\Â=\active\defÂ{\neg}        
\catcode`\µ=\active\defµ{\mu}         
\catcode`\¿=\active\def¿{{\o}}        
\catcode`\¹=\active\def¹{\pi}         
\catcode`\Ï=\active\defÏ{{\oe}}       
\catcode`\§=\active\def§{{\ss}}       
\catcode`\ =\active\def {\dagger}     
\catcode`\Ã=\active\defÃ{\sqrt}       
\catcode`\·=\active\def·{\Sigma}      
\catcode`\Å=\active\defÅ{\approx}     
\catcode`\½=\active\def½{\Omega}      
\catcode`\£=\active\def£{{\it\$}}     
\catcode`\°=\active\def°{\infty}      
\catcode`\¤=\active\def¤{{\S}}        
\catcode`\¦=\active\def¦{{\P}}        
\catcode`\¥=\active\def¥{\bullet}     
\catcode`\»=\active\def»{\leavevmode\raise.585ex\hbox{\b a}}      
\catcode`\¼=\active\def¼{\leavevmode\raise.6ex\hbox{\b o}}        
\catcode`\­=\active\def­{\not=}       
\catcode`\²=\active\def²{\leq}        
\catcode`\³=\active\def³{\geq}        
\catcode`\Ö=\active\defÖ{\div}        
\catcode`\É=\active\defÉ{{\dots}}     
\catcode`\¾=\active\def¾{{\ae}}       
\catcode`\Ç=\active\defÇ{\og}         
\catcode`\Ò=\active\defÒ{``}          
\catcode`\Á=\active\defÁ{!`}          
\catcode`\¢=\active\def¢{\rlap/c}     
\catcode`\Ô=\active\defÔ{`}           
\catcode`\Õ=\active\defÕ{'}           


\catcode`\=\active\def{{\AA}}       
\catcode`\'=\active\def'{\c C}        
\catcode`\¯=\active\def¯{{\O}}        
\catcode`\¸=\active\def¸{\Pi}         
\catcode`\Î=\active\defÎ{{\OE}}       
\catcode`\®=\active\def®{{\AE}}       
\catcode`\×=\active\def×{\diamond}    
\catcode`\¡=\active\def¡{\accent'27}  
\catcode`\Ó=\active\defÓ{''}          
\catcode`\±=\active\def±{\pm}         
\catcode`\È=\active\defÈ{\fg}         
\catcode`\À=\active\defÀ{?`}          
\catcode`\Ð=\active\defÐ{--}          
\catcode`\Ñ=\active\defÑ{---}         


\catcode`\Š=\active\defŠ{\"a}        
\catcode`\'=\active\def'{\"e}        
\catcode`\•=\active\def•{\"{\i}}     
\catcode`\š=\active\defš{\"o}        
\catcode`\Ÿ=\active\defŸ{\"u}        
\catcode`\Ø=\active\defØ{\"y}        
\catcode`\å=\active\defå{\^A}        
\catcode`\€=\active\def€{\"A}        
\catcode`\…=\active\def…{\"O}        
\catcode`\†=\active\def†{\"U}        
\catcode`\‡=\active\def‡{\'a}        
\catcode`\Ž=\active\defŽ{\'e}        
\catcode`\'=\active\def'{\'{\i}}     
\catcode`\—=\active\def—{\'o}        
\catcode`\œ=\active\defœ{\'u}        
\catcode`\ƒ=\active\defƒ{\'E}        
\catcode`\æ=\active\defæ{\^E}        
\catcode`\é=\active\defé{\`E}        %
\catcode`\ˆ=\active\defˆ{\`a}        
\catcode`\=\active\def{\`e}        
\catcode`\"=\active\def"{\`{\i}}     
\catcode`\˜=\active\def˜{\`o}        
\catcode`\=\active\def{\`u}        
\catcode`\Ë=\active\defË{\`A}        
\catcode`\‹=\active\def‹{\~a}        
\catcode`\–=\active\def–{\~n}        
\catcode`\›=\active\def›{\~o}        
\catcode`\Ì=\active\defÌ{\~A}        
\catcode`\"=\active\def"{\~N}        
\catcode`\Í=\active\defÍ{\~O}        
\catcode`\‰=\active\def‰{\^a}        
\catcode`\=\active\def{\^e}        
\catcode`\"=\active\def"{\^{\i}}     
\catcode`\™=\active\def™{\^o}        
\catcode`\ž=\active\defž{\^u}        

\let\optionkeymacros\null

\ifx\montrerlabels\oui
\input montrerlabels.tex
\fi
\def\paradouze{oui}
\anglais

\def\QQ{{\scal Q}}
\def\gR{{\got R}}
\def\sS{{\scal S}}
\hautspages{G. Tenenbaum}{On a family of arithmetic series related to the Mšbius function}\dateheure
\titrecentre{On a family of arithmetic series related to the Mšbius function}
\bigskip\medskip
\centerline{GŽrald Tenenbaum}
\bigskip\bigskip
{\leftskip75mm\it \obeylines
To George Andrews and Bruce Berndt,
as a friendly token of companionship
\par } 
\bigskip\bigskip
{\eightpoint\leftskip1cm\rightskip1cm
\noi{\bf Abstract.} Let $P^-(n)$ denote the smallest prime factor of a natural integer $n>1$. Furthermore let $\mu$ and $\omega$ denote respectively the Mšbius function and the number of distinct prime factors function. We show that, given any set $\P$ of prime numbers with a natural density, we have $\sum_{P^-(n)\in\P}\mu(n)\omega(n)/n=0$ and provide a effective estimate for the rate of convergence. This extends a recent result of Alladi and Johnson, who considered the case when $\P$ is an arithmetic progression.
 \PAR
\medskip
\noi
{\bf Keywords:} Mšbius function, Mšbius inversion, Perron's formula, saddle-point estimates, contour integration.\par
\smallskip 
\noi \bf 2020 Mathematics Subject Classification: \rm primary   11N37; secondary 11N25, 11N56.\par }
\bigskip
\medskip
\paraun{Introduction and statements}
Let $P^-(n)$ (resp. $P^+(n)$) denote the smallest (resp. the largest) prime factor of a natural integer $n>1$ and put $P^-(1):=\infty$ (resp. $P^+(1):=1)$. Furthermore, let $\mu$ and $\omega$ denote respectively the Mšbius function and the number of distinct prime factors function. \par 
In a recent paper \citer{AJ24}, Alladi and Johnson proved that, for given integers $k$, $\ell$, such that $(k,\ell)=1$, we have
$$\sum_{\di{n\leqslant x}{P^-(n)\md \ell k}}{\mu(n)\omega(n)\over n}\ll{(\log_2x)^{5/2}\over \sqrt{\log x}}\qquad (x\geqslant 3),$$ and consequently
that
$$\sum_{P^-(n)\md \ell k}{\mu(n)\omega(n)\over n}=0.\eqdef{AJser}$$
Their proof rests significantly on the \pnt\ for arithmetic progressions and on a duality identity due to Alladi \citer{Al77}, connecting small and large prime factors via Mšbius inversion. The purpose of this note is to investigate to what extent \eqref{AJser} depends on the subset of the primes appearing in the summation condition. We obtain the following result.  Here  and in the sequel we use the notation $u:=(\log x)/\log y$ $(x\geqslant y\geqslant 2)$, and we let $\log_k$ denote the $k$-fold iterated logarithm.
\vskip-3mm
\Propt{thgen}{Let $\P$ be a set of prime numbers satisfying, for suitable $\delta\in[0,1]$, 
$$\varepsilon_\P(t):={1\over t}\sum_{\di{p\leqslant t}{p\in\P}}\log p-\delta =o(1)\qquad  (t\to\infty).\eqdef{dens}$$
Then 
$$\sum_{P^-(n)\in\P}{\mu(n)\omega(n)\over n}=0.\eqdef{ser}$$
Moreover, for any fixed $b>5/3$ and uniformly for $\e^{(\log_2x)^b}\leqslant y\leqslant \sqrt{x}$, we have
$$\sum_{\di{n\leqslant x}{P^-(n)\in\P}}{\mu(n)\omega(n)\over n}\ll\varepsilon_\P^*(y)\log u+{1\over u},\eqdef{eff}$$
where $\varepsilon_\P^*(y):=\sup_{t>y}|\varepsilon_\P(t)|$.}
\par\goodbreak \medskip
\rem Quasi-optimal choices for $y$ yield that the upper bound in \eqref{eff} is, with arbitrary constants $\sigma>0$, $0<\tau<3/5$,
$$\ll\cases{ \dsp{\log_3x\over (\log_2x)^{\sigma}}& if $\varepsilon_\P^*(y)\ll 1/(\log_2 y)^{\sigma}$\cr
\dsp {(\log_2x)^{1/(1+\sigma)}\over (\log x)^{\sigma/(1+\sigma)}} & if $\varepsilon_\P^*(y)\ll 1/(\log y)^{\sigma}$,  \cr
\dsp{(\log_2x)^{1/\tau}\over \log x} & if $\varepsilon_\P^*(y)\ll \e^{-(\log y)^{\tau}}$.
}$$
The last case  covers that of an arithmetic progression.
\medskip\goodbreak
Let $\PP$ denote the set of all prime numbers. We note that \eqref{ser} does not hold for an arbitrary set of primes.  As suggested by Alladi in private communication, selecting
$$\P:=\PP\cap\cup_{j\geqslant 1}\big]\sqrt{x_j},x_j]\eqdef{Pcontrex}$$
for sufficiently rapidly increasing sequence $\{x_j\}_{j=1}^{\infty}$ furnishes a counter-example. Indeed, a straightforward consequence of \eqref{P0} {\it infra} is that
$$\liminf_{x\to\infty}\sum_{\di{n\leqslant x}{P^-(n)\in\P}}{\mu(n)\omega(n)\over n}\leqslant -\log 2.\eqdef{contrex}$$
See Remark 2.3 below.
\medskip\bigskip
\paraun{Proof of \ref{thgen}}
Let $y\in[2,x]$ be a parameter at our disposal, and put $\P_y:=\P\cap[2,y]$. We first aim at estimating the contribution from $\P_y$ to the sum~\eqref{eff} when~$y$ is sufficiently small in front of~$x$.The following lemma will be useful.
We define
$$\chi_{\P}(n):=\1_{\P}\big(P^-(n)\big),\quad g_y(n):=\sum_{m|n}\chi_{\P_y}(m)\mu(m)\omega(m)\quad \big(n\geqslant 1,\,y\geqslant 2\big),\eqdef{defgy}$$
 let $\sS$ denote the set of prime powers, and write $\sS^*:=\sS\cup\{1\}$
\Propl{lgy}{We have
$$g_y(n)=-\1_{\P_y}\big(P^+(n)\big)+\sum_{\tri{rs=n}{P^+(r)\in\P_y}{P^+(r)<P^-(s)=P^+(s)}}1\quad(n\geqslant 1).\eqdef{fgy}$$
In particular, for all $n\geqslant 1$, we have $|g_y(n)|\leqslant 1$ and
$$|g_y(n)|\leqslant \sum_{\tri{rs=n}{P^+(r)\leqslant y,\,s\in\sS^*}{P^-(s)>P^+(r)}}1.\eqdef{majgy}$$}
\nid
Put $\chi(n,y):=\1_{[1,y]}\big(P^+(n)\big)\ (n\geqslant 1).$ 
  Using the representation $n=ab$ with $\chi(a,y)=1$, $P^-(b)>y$, we have
$$\eqalign{g_y(n)&=\sum_{d|a,\,t|b}\chi_{\P_y}(d)\mu(d)\mu(t)\{\omega(d)+\omega(t)\}\cr
&=\sum_{d|a}\chi_{\P_y}(d)\mu(d)\omega(d)\sum_{t|b}\mu(t)+\sum_{d|a}\chi_{\P_y}(d)\mu(d)\sum_{t|b}\mu(t)\omega(t).\cr}$$
However
$$\eqalign{&\sum_{t|b}\mu(t)=\chi(n,y),\quad\sum_{t|m}\mu(t)\omega(t)=\Big[{\dd(1-z)^{\omega(m)}\over \dd z}\Big]_{z=1}=-\1_{\sS(m)}\quad(m\geqslant 1),\cr
&\sum_{d|a}\chi_{\P_y}(d)\mu(d)\omega(d)=-\sum_{p|a}\1_{\P_y}(p)\sum_{\di{d|a/p}{P^-(d)>p}}\mu(d)\{1+\omega(d)\}\cr
&\hskip30mm=-\1_{\P_y}\big(P^+(a)\big)+\sum_{\tri{a=rm}{P^+(r)\in\P_y}{P^+(r)<P^-(m)=P^+(m)\leqslant y}}1,
\cr
&\sum_{d|a}\chi_{\P_y}(d)\mu(d)=-\sum_{\di{p\in\P_y}{p|a}}\sum_{\di{d|a/p}{P^-(d)>p}}\mu(d)=-\1_{\P_y}\big(P^+(a)\big).\cr
}$$
This plainly implies \eqref{fgy}.
\qed
\medskip
We are now in a position to estimate the quantity
$$A^-(x,y):=\sum_{\di{n\leqslant x}{P^-(n)\in\P_y}}{\mu(n)\omega(n)\over n}\cdot$$
\Propl{lA-}{Without any hypothesis on $\P$ and uniformly for $2\leqslant y\leqslant \sqrt{x}$, we have
$$A^-(x,y)\ll{1\over u}\cdot\eqdef{P0}$$}
\nid Interpreting definition \eqref{defgy} as the Dirichlet convolution $g_y=\chi_{\P_y}\mu\omega*\1$, we get $\chi_{\P_y}\mu\omega=g_y*\mu$ by M\"obius inversion. Appealing to a strong form of the \pnt\ (see \eg\ \citeplus{BT04}{th. 8.17} or \citeplus{Te15}{ex. 178}), it follows that
$$\eqalign{A^-(x,y)&=\sum_{n\leqslant x}{\chi_{\P_y}(n)\mu(n)\omega(n)\over n}\cr
&=\sum_{d\leqslant x}{g_y(d)\over d}\sum_{m\leqslant x/d}{\mu(m)\over m}\ll\sum_{d\leqslant x}{|g_y(d)|\over d}\e^{-c\sqrt{\log x/d}}\cr}$$
for a suitable absolute constant $c>0$.
Now by \eqref{majgy} we have, for $2\leqslant y\leqslant D^{2/3}$,
$$\sum_{d\leqslant D}|g_y(d)|\ll \sum_{\di{r\leqslant D}{P^+(r)\leqslant y}}1+\sum_{\di{rP^+(r)\leqslant D}{P^+(r)\leqslant y}}{D\over r\log (DP^+(r)/r)}\ll {D\log y\over \log D},\eqdef{majsomgy}$$
where the $r$-sums have been estimated using the estimate \citeplus{Te15}{th. III.5.1}
$$\sum_{\di{r\leqslant t}{P^+(r)\leqslant y}}1\ll t\e^{-(\log t)/2\log y}\quad(t\geqslant 1,\,y\geqslant 2).$$
Recalling notation $u:=(\log x)/\log y$, we get
$$\sum_{d\leqslant x}{|g_y(d)|\over d}\e^{-c\sqrt{\log x/d}}\ll {1\over u}\qquad \big(2\leqslant y\leqslant \sqrt{x}\big),$$
where the estimate is obtained by splitting the summation at $x^{3/4}$, appealing to the inequality $|g_y(d)|\leqslant 1$ for the lower range and performing partial summation resting on~\eqref{majsomgy} for the upper range.\par
This implies \eqref{P0}, as required.
\qed \goodbreak
\medskip
\noi {\bf Remark 2.3}. To prove \eqref{contrex}, observe that if $\P$ is defined by \eqref{Pcontrex} then
$$A^-\big(x_j,\sqrt{x_j}\big)\ll{\log x_{j-1}\over \log x_j}$$
by \ref{lA-} and
$$\sum_{\tri{n²x_j}{P^-(n)\in\P}{P^-(n)>\sqrt{x_j}}}{\mu(n)\omega(n)\over n}=\sum_{\sqrt{x_j}<p\leqslant x_j}{-1\over p}=-\log 2+o(1)\quad(j\to\infty).$$
This implies \eqref{contrex} on selecting $\{x_j\}_{j=1}^{\infty}$ tending to infinity sufficiently fast.
\avance{\enonce}
\medskip
To complete the proof of \ref{thgen}, it remains to estimate the contribution from $\QQ_y:=\P\sset\P_y$ to the sum \eqref{eff}, viz.
$$A^+(x,y):=\sum_{\di{n\leqslant x}{P^-(n)\in\QQ_y}}{\mu(n)\omega(n)\over n}\cdot$$ 
This is the purpose of the following statement. Here and throughout, we let $\gamma$ denote Euler's constant.
\Propl{estA+}{Let $b>5/3$. Unformly for $x\geqslant 3$, $\e^{(\log_2x)^b}\leqslant y\leqslant \sqrt{x}$, we have
$$A^+(x,y)= {\delta\e^{\gamma}\over u}+ O\bigg(\varepsilon_\P^*(y)\log u+{1\over u^{9/5}}+{u^{6/5}\over \log x}\bigg),\eqdef{estA+1}$$
and
$$A^+(x,y)={\delta\e^{\gamma}\over u}+O\bigg(\varepsilon_\P^*(y)\log_2x+{1\over u^{9/5}}+{1\over \log x}\bigg).\eqdef{estA+2}$$}
\nid Both estimates will be proved by interpreting $A^+(x,y)$ as the derivative at $z=1$ of the polynomial
$$A(x,y;z):=\sum_{\di{n\leqslant x}{P^-(n)\in\QQ_y}}{\mu(n)z^{\omega(n)}\over n}$$
and estimating this quantity by the Selberg-Delange method.
\par 
We first consider \eqref{estA+1}, which turns out to be the more delicate of the two.
\par 
 Let the letters $p$ and~$q$ denote prime numbers.
For  $|z-1|\leqslant 1/5$, $w\geqslant 1$, $\Re s>1$, define
$$\eqalign{
G(s;w,z)&:=\prod_{q>w}\Big(1-{1\over q^s}\Big)^{-z}\Big(1-{z\over q^s}\Big),\cr
F(s;w,z)&:=\sum_{P^-(n)>w}{z^{\omega(n)}\mu(n)\over n^s}=\prod_{q>w}\Big(1-{z\over q^s}\Big)=\prod_{q\leqslant w}\Big(1-{1\over q^s}\Big)^{-z}{G(s;w,z)\over \zeta(s)^{z}}\cdot\cr}$$
Then
$$\HH(s;y,z):=\sum_{P^-(n)\in\QQ_y}{z^{\omega(n)}\mu(n)\over n^s}=-z\sum_{p\in\QQ_y}{F(s;p,z)\over p^s}\cdot\eqdef{defHsyz}$$
\par 
\par 
By a variant of Perron's formula \citeplus{Te15}{lemma II.2.6}, there exist two constants $\alpha$ and $\beta$ such that, writing
$$k(s):={1\over s}+{\alpha\over s+1}+{\beta\over s+2}\big(s\in\CC\sset\{-2,-1,0\}\big),\quad g(t):=\1_{[1,\infty]}(t)\Big\{1+{\alpha\over t}+{\beta\over t^2}\Big\}\quad(t>0),$$ we have, uniformly for $v>0$, $\kappa>0$,
$${1\over 2\pi i}\int_{\kappa-i}^{\kappa+i}k(s)v^s\d s=g(v)+O\Big({v^\kappa\over 1+(\log v)^2}+{\kappa v^\kappa}\Big).$$
\goodbreak
We infer that, for $|z|=r$,  $\kappa:=1/\log x$,
$$A(x,y;z):=\sum_{\di{n\leqslant x}{P^-(n)\in\QQ_y}}{\mu(n)z^{\omega(n)}\over n}={1\over 2\pi i}\int_{\kappa-i}^{\kappa+i}\HH(s+1,y;z)k(s)x^s\d s+O\bigg(\sum_{1\leqslant j\leqslant 4}R_j\bigg),$$
with
$$\eqalign{R_1&:=\sum_{P^-(n)\in\QQ_y}{\mu(n)^2r^{\omega(n)}\over n^{\kappa+1} \{1+\log (x/n)^2\}},\quad R_2:=\kappa\sum_{P^-(n)\in\QQ_y}{\mu(n)^2r^{\omega(n)}\over n^{\kappa+1}},\cr
R_3&:={1\over x}\sum_{\di{n\leqslant x}{P^-(n)\in\QQ_y}}\mu(n)z^{\omega(n)},\hskip17mm R_4:={1\over x^2}\sum_{\di{n\leqslant x}{P^-(n)\in\QQ_y}}n\mu(n)z^{\omega(n)}.\cr}$$
\par
 We readily have $R_2\ll u^r/\log x$.
 To estimate $R_1$,  first consider the contribution, say~$R_{11}$, of those integers $n$ 
such that $|\log (x/n)|>1$. Summing over dyadic intervals and appealing to standard bounds for averages of non-negative arithmetic functions, \eg\ \citeplus{Te15}{th. III.3.5}, we see that \hbox{$R_{11}\ll u^r/\log x$}.  The complementary contribution $R_{12}$ is a sum over $[x/\e,\e x]$. It is readily evaluated by  applying the same standard bounds. This yields again $R_{12}\ll u^r/\log x$.
The terms $R_3$ and $R_4$ may be estimated trivially by bounding $\mu(n)$ by $\mu(n)^2$ and $z^{\omega(n)}$ by $r^{\omega(n)}$. This still furnishes $R_3+R_4\ll u^r/\log x$. \par \goodbreak
We may finally state that
$$A(x,y;z)={1\over 2\pi i}\int_{\kappa-i}^{\kappa+i}\HH(s+1,y;z)k(s) x^s\d s+O\Big({u^r\over \log x}\Big).\eqdef{Perron}$$
\par 
Define $$J(s):=\int_0^\infty \e^{-s-t}{\dd t\over s+t}\quad\big(\Re s>0\big),\qquad L_\varepsilon(t):=\e^{(\log t)^{3/5-\varepsilon}}\quad(\varepsilon>0,\,t\geqslant 2).$$ When $\Re s>0$, $p>y$, $s_p:=s\log p$,  \citeplus{Te15}{lemma III.5.16} yields, for any fixed $\varepsilon>0$,
 $$F(s+1;p,z)= \e^{-zJ(s_p)}\Big\{1+O\Big({1\over L_\varepsilon(y)}\Big)\Big\}\qquad \Big(|\Im s|\leqslant L_\varepsilon(y)\Big).\eqdef{evalF}$$  
 (This is proved using the Korobov-Vinogradov zero free region and accounts for the exponent $3/5$ in the definition of $L_\varepsilon(y)$.) 
\par 
Let us insert \eqref{evalF} into \eqref{defHsyz} and then into \eqref{Perron} keeping in mind the hypothesis $\log y\geqslant (\log_2x)^b$. Using the estimate $|\e^{-J(s)}|\asymp \min(|s|,1)$ $(\Re s\geqslant -1)$ proved in \hbox{\citeplus{HT93}{lemma~2}} in the form
$$|\e^{-zJ(s)}|\asymp \min(|s|^{\Re z},1)\qquad (\Re s\geqslant -1,\,\Re z>0),\eqdef{ezJ}$$
we obtain
$$A(x,y;z)={1\over 2\pi i}\int_{\kappa-i}^{\kappa+i}B(s;y,z)k(s)x^s\d s+O\Big({u^{r}\over \log x}\Big),\eqdef{P2}$$
with
$$B(s;y,z):=\sum_{p\in\QQ_y}{-z\e^{-zJ(s_p)}\over p^{s+1}}\cdot$$
\par 
Now we introduce the remainder $$R(t):=t\varepsilon_\P(t)=\sum_{\di{p\leqslant t}{p\in \P}}\log p-\delta t=o(t)\qquad (t>1).$$
Taking into account that 
$J'(s) =-\e^{-s}/s$, we get
$$B(s;y,z)=D(s;y,z)-z\int_{y}^{\infty}{\e^{-zJ(s_t)}\over t^{s+1}\log t}\d R(t),\eqdef{evalB}$$
with
$$\eqalign{D(s;y,z)&:=-\delta z\int_{y}^{\infty}{\e^{-zJ(s_t)}\over t^{s+1}\log t}\d t=\delta\int_{s_y}^{s\infty}{-z \e^{-zJ(v)}\over v\e^v}\d v=\delta\big\{1-\e^{-zJ(s_y)}\big\}.\cr}$$
\par 
Carrying back into \eqref{P2}, we obtain
$$\eqalign{A(x,y;z)&={\delta\over 2\pi i}\int_{\kappa-i}^{\kappa+i}\big\{1-\e^{-zJ(s_y)}\big\}k(s)x^s\d s+O\Big(\gR_\P(x,y;z)+{u^{r}\over \log x}\Big),\cr}\eqdef{P3}$$
with $$\gR_\P(x,y;z):=\int_y^\infty{\lambda_x(t)\over t\log t}\d R(t),\quad \lambda_x(t):=\int_{\kappa-i}^{\kappa+i}{\e^{-zJ(s_t)}}\Big({x\over t}\Big)^sk(s)\d s.$$
\par \goodbreak
By \eqref{ezJ}, we have, for $t>y$,
$$\eqalign{\lambda_x(t)&\ll\int_{\kappa-i}^{\kappa+i}\Big({x\over t}\Big)^s\min(|s|\log t,1)^{\Re z}{|\dd s|\over |s|}\ll \Big({x\over t}\Big)^\kappa\log_2t ,\cr
\lambda_x'(t)&=\int_{\kappa-i}^{\kappa+i}x^sk(s){\dd\over \dd t}\Big({\e^{-zJ(s_t)}\over t^s}\Big)\d s=\int_{\kappa-i}^{\kappa+i}{x^ssk(s)\e^{-zJ(s_t)}\over t^{s+1}}\Big\{-1+{z\over st^s\log t}\Big\}\d s\cr&\ll{x^\kappa\over t^{\kappa+1}}\bigg\{1+\int_{\kappa-i}^{\kappa+i}{\min(|s|\log t,1)^{\Re z}\over |s|\log t}|\dd s|\bigg\}\ll{x^\kappa\over t^{\kappa+1}}\cdot\cr}$$
Partial integration hence furnishes
$$\gR_\P(x,y;z)\ll\varepsilon_\P^*(y)\log u.\eqdef{evalRP}$$
Now we know \citeplus{Te15}{(III.5.41)} that $\e^{-J(s)}=s\wh\varrho(s)$, where
$$\wh\varrho(s):=\int_0^\infty \varrho(v)\e^{-sv}\d v,$$
an entire function, is the Laplace transform of the Dickman function  $\varrho(v)$. Therefore, assuming with no loss of generality that $x\in\dm+\N$, the main term, say $M$, in \eqref{P3} satisfies
$$\eqalign{M&:=\delta-{\delta\over 2\pi i}\int_{\kappa-i}^{\kappa+i}\{s\log y\}^{z}\wh\varrho(s\log y)^zk(s)x^s\d s+O\Big({1\over \log x}\Big)\cr
&=\delta-{\delta\over 2\pi i}\int_{1/u-i\log y}^{1/u+i\log y}w^{z-1}\wh\varrho(w)^zk_y(w)\e^{u w}\d w+O\Big({1\over \log x}\Big),\cr}\eqdef{M1}$$
where we have put
$$k_y(w):=1+{\alpha w\over w+\log y}+{\beta w\over w+2\log y}\cdot$$
The last integral may be evaluated on replacing the integration segment by a truncated Hankel contour around $\big[-\dm,1/u\big]$,\note{That is the path consisting of the circle $|w|=1/u$ excluding the point $w=-1/u$ and the segment $[-1/2,-1/u]$ covered twice with respective arguments $+\pi$ and $-\pi$.} concatenated with two vertical segments \hbox{$\big[-\dm,-\dm\pm i\log y\big]$} and two horizontal segments $\big[-\dm\pm i\log y,1/u\pm i\log y\big]$. Appealing for instance to \hbox{\citeplus{Te15}{cor. II.0.18}} to take care of the truncation, and noting that $\wh\varrho(w)=\e^\gamma+O(w)$ for $w\ll1$, we see that the contribution of the Hankel contour is $${\e^{\gamma z}\over \Gamma(1-z)u^z}+O\Big({1\over u^{\Re z+1}}\Big).$$
\par \goodbreak
On the vertical parts of the contour, we have $w^z\wh\varrho(w)^z\ll1$ if $|\Im w|\leqslant 1$ and $w^z\wh\varrho(w)^z=1+O(1/w)$ if $|\Im w|>1$. The first range contributes $\ll\e^{-u/2}$ to the integral. The contribution of the second may be estimated on noting that, writing $I:=[-\log y,-1]\cup[1,\log y]$, we have
$$\int_{I}{\e^{-u/2+i\tau u}\over -1/2+i\tau}\d\tau=\int_1^{\log y}{-\e^{-u/2}\{\cos(\tau u)+2\tau\sin(\tau u)\}\over 1/4+\tau^2}\d\tau\ll \e^{-u/2}$$
by the second mean-value theorem. Finally, the contribution of the horizontal parts is trivially 
$$\ll \int_{-\1/2}^{1/u}{\e^{\sigma u}\over \log y}\d\sigma\ll{1\over \log x}\cdot$$
Thus 
$$M=\delta-{\delta\e^{\gamma z}\over \Gamma(1-z)u^z}+O\Big({1\over u^{\Re z+1}}+{1\over \log x}\Big).$$
Gathering our estimates, we arrive at
$$A(x,y;z)=\delta -{\delta\e^{\gamma z}\over \Gamma(1-z)u^z}+O\Big(\varepsilon_\P^*(y)\log u+{1\over u^{\Re z+1}}+{u^r\over \log x}\Big).$$
Differentiating at $z=1$  using Cauchy's integral formula, we get \eqref{estA+1} as required since $\Re z\geqslant 4/5$, $ r\leqslant 6/5$.
\par \smallskip
We now turn our attention to proving \eqref{estA+2}. 
To this end, we appeal to a standard Perron formula  \citeplus{Te15}{th. II.2.3}, viz.
$$A(x,y;z)={1\over 2\pi i}\int_{\kappa-iT}^{\kappa+iT}\HH(s+1,y;z){x^s\over s}\d s+O\bigg(\sum_{P^-(n)\in\QQ_y}{\mu(n)^2r^{\omega(n)}\over n^{\kappa+1} (1+T|\log (x/n|)}\bigg),$$
where $r=|z|\in[4/5,6/5]$, $|z-1|\leqslant 1/5$.
 Those integers $n$ such that $|\log (x/n)|>1$ contribute $\ll u^{r}/(T\log x)\ll(\log x)^{r-1}/T$ to the error term. Splitting the summation range of the complementary sum into intervals $]x+hx/T, x+(h+1)x/T]$ $(|h|\leqslant T)$ and applying Shiu's theorem \citer{Sh80} for short sums of multiplicative functions, we obtain that it is $\ll (u^r\log T)/(T\log x)\ll (\log x)^{r-1}(\log T)/T$ provided, say, $2\leqslant T\leqslant \sqrt{x}$.  
 \par  Select $T:=(\log x)^{r+1}$, so that, in view of hypothesis $\log y\geqslant (\log_2x)^b$, we have $T \leqslant L_\varepsilon(y)$ for suitable $\varepsilon>0$. This yields
$$\eqalign{A(x,y;z)&={1\over 2\pi i}\int_{\kappa-iT}^{\kappa+iT}\HH(s+1,y;z){x^s\over s}\d s+O\Big({(\log x)^{r}\over T}\Big)\cr
&={1\over 2\pi i}\int_{\kappa-iT}^{\kappa+iT}B(s;y,z){x^s\over s}\d s+O\Big({1\over \log x}\Big),\cr}\eqdef{P2+}$$
by \eqref{evalF}. \par \goodbreak
From \eqref{evalB}, we get
$$A(x,y;z)={\delta\over 2\pi i}\int_{\kappa-iT}^{\kappa+iT}\big\{1-\e^{-zJ(s_y)}\big\}{x^s\over s}\d s+O\Big(\gR_\P^+(x,y;z)+{1\over \log x}\Big),\eqdef{P3+}$$
with
$$\gR_\P^+(x,y;z):=\int_y^\infty{\nu_x(t)\over t\log t}\d R(t),\quad \nu_x(t):=\int_{\kappa-iT}^{\kappa+iT}{\e^{-zJ(s_t)}}\Big({x\over t}\Big)^s{\dd s\over s}\cdot$$
Appealing to the estimates $$\big|\e^{-zJ(s)}\big|\asymp \min\big(|s|^{\Re z},1\big),\quad\e^{-zJ(s)}=1+O(1/s)\quad(\Re z>0,\, \Re s>0),$$ we get, for $t>y$, keeping in mind  the hypothesis $\log y\geqslant (\log_2x)^b$,
$$\eqalign{\nu_x(t)&\ll\Big({x\over t}\Big)^\kappa \log_2x,\cr 
\nu_x'(t)&=\int_{\kappa-iT}^{\kappa+iT}\e^{-zJ(s_t)}{x^s\over t^{s+1}}\Big\{{z\over t^ss\log t}-1\Big\}\d s\cr
&=-\int_{\kappa-iT}^{\kappa+iT}\e^{-zJ(s_t)}{x^s\over t^{s+1}}\d s+O\bigg({x^\kappa\log_2x\over t^{\kappa+1}\log y}\bigg)\cr
&=-\int_{\kappa-iT}^{\kappa+iT}{x^s\over t^{s+1}}\d s+O\bigg(\int_{\kappa-iT}^{\kappa+iT}\vbs{x^s\over t^{s+1}s\log t}|\dd s|+{x^\kappa\over t^{\kappa+1}}\bigg)\cr
&\ll {x^\kappa T\over t^{\kappa+1}(1+T|\log (x/t)|)}+{x^\kappa\over t^{\kappa+1}},\cr}$$
where we used the bound $$\int_{\kappa-iT}^{\kappa+iT}w^s\d s\ll {w^\kappa T\over 1+T|\log w|}\qquad (w>0,\,\kappa>0,\,T\geqslant 1).$$
From the above estimates, we infer  that
$\gR_\P^+(x,y;z)\ll \varepsilon_\P^*(y)\log_2x.$
\par \smallskip
Assuming as before that $x\in\dm+\N$, the main term $M^+$ in \eqref{P3+} satisfies
$$\eqalign{M^+&=\delta-{\delta\over 2\pi i}\int_{1/u-iT\log y}^{1/u+iT\log y}w^{z-1}\wh\varrho(w)^z\e^{u w}\d w+O\Big({1\over \log x}\Big)\cr
&=\delta-{\delta\e^{\gamma z}\over\Gamma(1-z)u^z}+O\Big({1\over u^{\Re z+1}}+{1\over \log x}\Big),\cr}$$
after deforming the integration segment and exploiting the relevant Hankel contour.
\par \goodbreak Finally, we may state that, for $|z-1|\leqslant 1/5$, we have
$$A(x,y;z)=\delta -{\delta\e^{\gamma z}\over \Gamma(1-z)u^z}+O\Big(\varepsilon_\P^*(y)\log_2x+{1\over u^{9/5}}+{1\over \log x}\Big).$$
Differentiating the above formula at $z=1$ furnishes \eqref{estA+2}.
\qed
\medskip
We are now  able to complete the proof of the effective estimate \eqref{eff}. This amounts to showing, that, in the stated range for $y$, we have
$$
A^-(x,y)+A^+(x,y)\ll\varepsilon_{\P}^*(y)\log u+1/u. \eqdef{estA}
$$
If $\log y\gg (\log x)^{6/11}$, then $u^{6/5}/\log x=(\log x)^{6/5}/\{u(\log y)^{11/5}\}\ll 1/u $, so \eqref{P0} and \eqref{estA+1} imply \eqref{estA}. If $\log y\ll (\log x)^{6/11}$, then $u\gg(\log x)^{5/11}$, so $\log_2x\ll\log u$, and \eqref{estA} follows from \eqref{P0} and \eqref{estA+2}.
\goodbreak\medskip\medskip
\paraun{Special cases}
We provide asymptotic formulae when $\P$ is either the set of all primes or a singleton. In these special cases, direct computations can be performed via  standard applications of the Selberg-Delange method, furnishing improved estimates.\note{This is due in particular to the fact that estimates like \eqref{evalF}, requiring information on the zeros of the zeta function, may be avoided.} Consequently we only sketch the main lines.\Propp{PV1}{We have
$$V_1(x):=\sum_{n\leqslant x}{\mu(n)\omega(n)\over n}\sim {-1\over \log x}\qquad (x\to\infty).\eqdef{V1}$$
}
\nid  Observe that, for $z\in\CC$, $|z|\leqslant \ft32$,
$$F_1(s,z):=\sum_{n\geqslant 1}{\mu(n)z^{\omega(n)}\over n^s}=\prod_{q}\Big(1-{z\over q^s}\Big)={G_1(s,z)\over \zeta(s)^z},$$
with
$$G_1(s,z):=\prod_{q}\Big(1-{z\over q^s}\Big)\Big(1-{1\over q^s}\Big)^{-z}.$$
Hence
$$V_1(x;z):=\sum_{n\leqslant x}{z^{\omega(n)}\mu(n)\over n}={1\over 2\pi i}\int_{\kappa-i\infty}^{\kappa+i\infty}{G_1(s+1,z)\over \{s\zeta(s+1)\}^{z}}{x^s\over s^{1-z}}\d s.$$
The main contribution arises from   a Hankel contour around $[-c,0]$ for arbitrary constant $c>0$. By Hankel's formula, we get, as $x\to\infty$,
$$V_1(x;z)\sim {G_1(1,z)\over \Gamma(1-z)(\log x)^{z}}={(1-z)G_1(1,z)\over \Gamma(2-z)(\log x)^z}\cdot$$
Hence
$$V_1(x)=\Big[{\dd V_1(x;z)\over \dd z}\Big]_{z=1}\sim{-G_1(1,1)\over \log x}={-1\over \log x},$$
as wanted.\qed
\goodbreak
\bigskip
Next consider the case of $\P$ being reduced to a single element. Write $$\zeta(s,y):=\prod_{q\leqslant y}\Big(1-{1\over q^s}\Big)^{-1}\qquad \big(\Re s>0,\,y\geqslant 2\big).$$
\Propp{PVp}{Let $p\in\PP$. We have
$$V_p(x):=\sum_{\di{n\leqslant x}{P^-(n)=p}}{\mu(n)\omega(n)\over n}\sim{\zeta(1,p)\over p\log x}\qquad (x\to\infty).\eqdef{Vp}$$}
\nid
Consider
$$F_p(s,z):=\sum_{\di{n\geqslant 1}{P^-(n)=p}}{\mu(n)z^{\omega(n)}\over n^s}={-z\over p^s}\prod_{q>p}\Big(1-{z\over q^s}\Big)={-zG_p(s,z)\over p^s\zeta(s)^{z}},$$
with now
$$G_p(s,z):=\prod_{q\leqslant p}\Big(1-{z\over q^s}\Big)^{-1}G_1(s,z).$$
\goodbreak
It follows that, as $x\to\infty$,
$$\eqalign{V_p(x;z):&=\sum_{\di{n\leqslant x}{P^-(n)=p}}{\mu(n)z^{\omega(n)}\over n}={-z\over 2\pi i}\int_{\kappa-i\infty}^{\kappa+i\infty}{G_p(s+1,z)\over \{s\zeta(s+1)\}^z}{x^s\over p^ss^{1-z}}\d s\cr&\sim{-z(1-z)G_p(1,z)\over \Gamma(2-z)(\log x/p)^z}\cdot\cr}$$
Differentiating at $z=1$ taking the zero of the numerator into account, we obtain \eqref{Vp}.
\qed\goodbreak
\bigskip
From the two propositions above, it follows that one cannot heuristically reconstruct~\eqref{V1} from \eqref{Vp}. This phenomenon is similar to that arising from the formulae
$$\eqalign{&\sum_{n\geqslant 1}{\mu(n)\log n\over n}=-1,
\quad\sum_{\di{n\geqslant 1}{P^-(n)=p}}{\mu(n)\log n\over n}={\zeta(1,p)\over p}\cdot}\eqdef{serlog}$$
\bigskip
\noi{\bf Acknowledgements.} The author expresses his gratitude to the referee for judicious advices and suggestions. He also thanks RŽgis de la Bretche for interesting exchanges on this problem, and Krishna Alladi for sharing an early version of his paper with Johnson and further friendly discussions on this topic.
\bigskip
\centerline{\twelvebf References}\bigskip
{\leftskip5mm\rightskip5mm\eightpoint{
\bibtem{Al77} K. Alladi, Duality between prime factors and the prime number theorem for arithmetic progressions, {\it J. Number Theory \bf 9} (1977), 436-451.\par 
\bibtem{AJ24} K. Alladi \& J. Johnson, Duality between prime factors and the prime number theorem for arithmetic progressions, II, preprint, July 2024.\par 
\bibtem{BT04} P.T. Bateman \& H. Diamond, Analytic number theory, volume 1 of Monographs in Number theory, An introductory course, Word Scientific, Hackensack, 2004.
\bibtem{HT93} A. Hildebrand \& G. Tenenbaum, On a class of difference differential equations arising in number
theory, {\it J. Anal. Math. \bf 61} (1993),  145--179.\par 
\bibtem{Sh80}
P. Shiu, A Brun-Titschmarsh theorem for multiplicative functions, {\it J. reine angew. Math.} {\bf 313} (1980), 161--170.
\par 
\bibtem{Te15} G. Tenenbaum, {\it Introduction to analytic and probabilistic number theory}, 3rd ed., Graduate Studies in Mathematics 163, Amer. Math. Soc. 2015.\par
\par }
}
\vskip 5mm
{\sevenrm\baselineskip9pt
G\'erald Tenenbaum\par
Institut \'Elie Cartan\par 
Universit\'e de Lorraine\par
 BP 70239\par
54506 Vand\oe uvre Cedex\par
 France
\smallskip
e-mail : \seventt gerald.tenenbaum@univ-lorraine.fr\par}

\bye